\documentclass[leqno,11pt,draft]{amsart}
\usepackage{amssymb}
\usepackage{amsmath}
\usepackage{enumerate}
\usepackage{amsfonts}

\newtheorem{teo}[equation]{Theorem}

\newtheorem{lema}[equation]{Lemma}

\numberwithin{equation}{section}

\newcommand{\R}{\mathbb R}
\newcommand{\Rie}{\mathcal R}
\newcommand{\Har}{H^1_L}
\newcommand{\Hardy}{H^1}
\newcommand{\LL}{L^1}
\newcommand{\M}{\mathcal{M}}
\newcommand{\8}{\infty}
\newcommand{\e}{\varepsilon}
\newcommand{\W}{\mathcal{W}}
\newcommand{\SSS}{\mathcal{S}}
\newcommand{\Q}{\mathcal{Q}}
\renewcommand{\d}{\partial}
\renewcommand{\dj}{\frac{\d}{\d x_j}}

\newcommand{\wt}{\widetilde}

\begin{document}

\title[Riesz transform characterization of Hardy spaces]
{Riesz transform characterization of Hardy\\
spaces associated with Schr\"odinger operators\\
with compactly supported potentials}

%???????????????????????????????????????????????????????????????????????????????
 \keywords{}
\begin{abstract} Let $L=-\Delta +V$ be a Schr\"odinger operator on $\R^d$,
$d\geq 3$. We assume that $V$ is a nonnegative, compactly
supported potential that belongs to $L^p(\R^d)$, for some
$p>d\slash 2$. Let $K_t$ be the semigroup generated by $-L$. We
say that an $L^1(\mathbb R^d)$-function $f$ belongs to the Hardy
space $\Har$ associated with $L$ if $\sup_{t>0} |K_t f|$ belongs
to $\LL(\R^d)$. We prove that $f\in \Har$ if and only if $R_j f
\in \LL(\R^d)$ for $j=1,...,d$, where $R_j= \dj L^{-1\slash 2}$
are the Riesz transforms associated with $L$.
\end{abstract}

%??????????????????????????????????????????????????????
\author[Jacek Dziuba\'nski ]{ Jacek Dziuba\'nski }
\address{Instytut Matematyczny\\
Uniwersytet Wroc\l awski\\
50-384 Wroc\l aw, Pl. Grunwaldzki 2/4\\
Poland} \email{jdziuban@math.uni.wroc.pl,
preisner@math.uni.wroc.pl}

\author[Marcin Preisner ]{ Marcin Preisner}

%????????????????????????????????????????????????//
\thanks{
 Supported by the European Commission Marie Curie Host Fellowship
for the Transfer of Knowledge "Harmonic Analysis, Nonlinear
Analysis and Probability" MTKD-CT-2004-013389 and by Polish
Government funds for science. }
%\date{\today}
\maketitle
%%%%%%%%%%%%%%%%%%%%%%%%%%%%%%%%%%%%%%%%%%%%%%%%%%%%%%%%%%%%%%%%%%%%%%%%%%%%%%%
\section{Introduction.}

 Let
 \begin{equation}\label{class_riesz}
  \mathcal R_j f(x)= c_d\lim_{\varepsilon \to 0}
 \int_{|x-y|>\varepsilon } \frac{x_j-y_j}{|x-y|^{d+1}} f(y)\,
 dy=\lim_{\varepsilon \to 0} \int_{\varepsilon
 }^{\varepsilon^{-1}}\int_{\mathbb R^d}
 \frac{\partial}{\partial x_j} P_t(x-y)f(y)\,
 dy\frac{dt}{\sqrt{t}},
 \end{equation}
 $j=1,2,...,d,$ be the classical Riesz transforms on $\mathbb R^d$. Here and
 subsequently
 $$ P_t(x-y)= (4\pi t)^{-d\slash 2}
 \exp \left(-|x-y|^2 \slash 4t\right)$$
 denotes the heat kernel. Clearly, for $f\in L^1(\mathbb R^d)$ the
 limits in (\ref{class_riesz}) exist  in the sense of
 distributions and define $\mathcal R_jf$ as a distribution. It was proved in
 Fefferman and Stein \cite{FS} (see also \cite{Stein}) that an $L^1(\mathbb R^d)$-function $f$ belongs
 to the classical Hardy space $H^1(\mathbb R^d)$ if and only if
 $\mathcal R_jf\in L^1(\mathbb R^d)$ for $j=1,...,d$. Moreover,
\begin{equation}\label{classical}
 \|f\|_{\LL(\mathbb R^d)} + \sum_{j=1}^d \|\Rie _j f\|_{\LL(\mathbb R^d)},
\end{equation}
 defines one of the possible norms in $H^1(\mathbb R^d)$.

In this paper we consider a Schr\"odinger operator $L=-\Delta +V$
on $\R^d$,  $d\geq 3$. We assume that $V$ is a nonnegative
function, $\mathrm{supp} \, V \subseteq B(0,1) =\{x\in \R: \,
|x|<1\}$,
 and  $V\in L^p(\R^d)$ for some $p>d\slash 2$.
Let $K_t$ be the semigroup generated by $-L$.
Since $V\geq 0$, by the Feynman-Kac formula, we have
\begin{equation}
\label{F-Kac} 0\leq K_t(x,y) \leq P_t(x-y)
\end{equation}
where $K_t(x,y)$ is the integral  kernel of the semigroup
$\{K_t\}_{t>0}$. Let
\begin{equation*}
\M f (x) = \sup_{t>0} |K_t f (x)|.
\end{equation*}
We say that an $L^1(\mathbb R^d)$-function $f$ belongs the Hardy
space $\Har$ if
\begin{equation*}
\|f\|_{\Har} = \|\M f\|_{\LL(\mathbb R^d)} < \8.
\end{equation*}

For $j=1,...,d$ let us define the Riesz transforms $R_j$ associated
with $L$ by setting
\begin{equation}\label{rieszdef}
R_j f = c'_d\frac{\d}{\d x_j} L^{-1\slash 2}f=\lim_{\varepsilon
\to 0}\int_{\varepsilon}^{1\slash \varepsilon}
\frac{\partial}{\partial x_j}K_t f\frac{dt}{\sqrt{t}},
\end{equation}
where the limit is understood in the sense of distributions. The
fact that for any $f\in L^1(\mathbb R^d)$ the operators
 \begin{equation}\label{R_j} R_j^{\varepsilon} f
 =\int_{\varepsilon}^{1\slash \varepsilon} \frac{\partial}{\partial x_j}
 K_tf\frac{dt}{\sqrt{t}}
 \end{equation}
 are well defined and the limit $\lim_{\varepsilon\to 0}  R^\varepsilon_j f$ exists
 in the sense of distributions will be discussed below.

  The
main result of this paper is the following.
\begin{teo}
\label{mainteo} Assume that  $f\in \LL(\mathbb R^d)$. Then  $f$ is
in the Hardy space $\Har$ if and only if $R_j f \in \LL(\mathbb
R^d)$ for every $j=1,...,d$. Moreover, there exists $C>0$ such
that
\begin{equation} \label{mainineq}
C^{-1} \|f\|_{\Har} \leq \|f\|_{\LL(\mathbb R^d)} + \sum_{j=1}^d
\|R_j f\|_{\LL(\mathbb R^d)} \leq C \|f\|_{\Har}.
\end{equation}
\end{teo}

 The Hardy spaces$H^1_L$ associated with the Schr\"odinger
 operators
 $L$ with compactly supported potentials were studied in \cite{DZ}.
 It was proved there that the elements of the space $H^1_L$ admit special
 atomic decompositions. Moreover, the space $H^1_L$ is isomorphic
 to the classical Hardy space $H^1(\mathbb R^d)$. To be more
 precise, let
 $$ \Gamma (x,y)=\int_0^\infty K_t(x,y)\, dt, \ \ \ \Gamma_0
 (x,y)=-\int_0^\infty P_t(x,y)\, dt,$$
 and denote
 $$ L^{-1} f(x) = \int_{\mathbb R^d} \Gamma (x,y)f(y)\, dy,\ \ \
 \Delta^{-1} f(x)= \int_{\mathbb R^d} \Gamma_0(x,y)f(y)\, dy.$$
The operators $(I-V\Delta^{-1})$ and $(I-VL^{-1})$ are bounded and
invertible on $\LL(\mathbb R^d)$, and
\begin{equation*}
I=(I-V\Delta^{-1})(I-VL^{-1})=(I-VL^{-1})(I-V\Delta^{-1}).
\end{equation*}
Moreover, $\left(I-VL^{-1}\right):\Har \to \Hardy(\mathbb R^d)$ is
an isomorphism (whose inverse is $(I-V\Delta^{-1})$) and
\begin{equation} \label{homeo}
\|\left(I-VL^{-1}\right) f \|_{\Hardy (\mathbb R^d)} \simeq
\|f\|_{\Har}
\end{equation}
for $f \in \Har$ (see \cite[Corollary 3.17]{DZ}).\\

 The proof of the special atomic decompositions presented in
 \cite{DZ} was based on the following identity
\begin{equation}\label{mainform}
 P_t(I-VL^{-1})=K_t-\int_0^t (P_t-P_{t-s})VK_s\, ds,-\int_t^\infty
P_tVK_s\, ds = K_t-W_t-Q_t,
\end{equation}
which comes from the
 perturbation formula
 $$ P_t=K_t+\int_0^\infty P_{t-s}VK_s\, ds.$$
 The formula (\ref{mainform}) will be also used here in the analysis of the
 integral (\ref{R_j}) for large $t$, while for $t$ small we shall
 use  its slightly different version, namely
\begin{equation}\label{pert}
P_t(I-VL^{-1})=K_t + \int_0^t P_{t-s} V K_s ds - P_t V L^{-1}=K_t
+ \widetilde{W_t}- \widetilde{Q_t}.
\end{equation}

In the case of Schr\"odinger operators with potentials $\mathcal
V\geq 0$, $\mathcal V  \not\equiv 0$, satisfying the reverse
H\"older inequality with the exponent $d\slash 2$ (which clearly
implies $\text{supp}\, \mathcal V=\mathbb R^d$ ),   Riesz transform
characterizations of the relevant Hardy spaces $H^1_{-\Delta
+\mathcal V}$ were obtained in \cite{DZ0}.

%%%%%%%%%%%%%%%%%%%%%%%%%%%%%%%%%%%%%%%%%%%%%%%%%%%%%%%%%%%%%%%%%%%%%%%%%%%%%%%
\section{Auxiliary estimates}
%Define
%$$\Q_j = \int_0^\8 \frac{\d}{\d x_j} Q_t \frac{dt}{\sqrt{t}}, \qquad
%\W_j = \int_0^\8 \frac{\d}{\d x_j} W_t \frac{dt}{\sqrt{t}}.$$
%In the following lemmas we will write $\Q,\W$ instead of $\Q_j, \W_j$ to simplify the notion, because all estimates will be true for all $j=1,...,d$. Finally, let us mention that $C$ will be used as a positive constant which may vary from line to line.
%In order to prove theorem \ref{mainteo} we will need some estimates for the operators $\Q$ and $\W$. The proofs of the following lemmas are very much in the spirit of \cite{DZ}.
In this section we will use notation  $f_t(x) = t^{-d\slash 2} f(\frac{x}{\sqrt{t}})$.

For $f \in \LL(\mathbb R^d)$ and $0<\varepsilon <1$ we define the
truncated Riesz transforms by setting
\begin{equation*}
R_j^{\e} f = \int_\e^{\e^{-1}} \dj K_t f \frac{dt}{\sqrt{t}},
\qquad \Rie_j^{\e} f = \int_\e^{\e^{-1}} \dj P_t f
\frac{dt}{\sqrt{t}}.
\end{equation*}
Denote
$$ G(x,y)=\int_0^\infty K_t(x,y)\frac{dt}{\sqrt{t}}, \ \ \
G_0(x,y)=\int_0^\infty P_t(x-y)\frac{dt}{\sqrt{t}}.$$
 Then
$G(x,y)\leq G_0(x,y)=c|x-y|^{-d+1}$ and, consequently, for
$\varphi$ from the Schwartz class $\mathcal S(\mathbb R^d)$ we
have
$$\lim_{\varepsilon \to 0} \langle R_j^\varepsilon f,\varphi
\rangle =-\int G(x,y)f(y)\frac{\partial}{\partial x_j} \varphi
(x)\, dy\, dx.$$ Hence $R_jf$ is a well defined distribution and
$$|\langle R_jf,\varphi\rangle | \leq C\| f\|_{L^1(\mathbb R^d)}\left(\left\|
\frac{\partial}{\partial x_j}\varphi\right\|_{L^1(\mathbb R^d)}+\left\|
\frac{\partial}{\partial x_j}\varphi\right\|_{L^\infty}\right).$$

 Using   \eqref{mainform} and \eqref{pert} we write
\begin{equation}\label{epsilon}
R_j^\e f = \Rie_j^\e (I-VL^{-1})f -\wt{\W_j^\e} f + \wt{\Q_j^\e}f
+ \W_j^\e f + \Q^\e_j f,
\end{equation}
where $\Q_j^\e$, $\wt{\Q_j^\e}$, $\W_j^\e$, and $\wt{\W_j^\e}$ are
operators with the following integral kernels
$$\Q_j^\e (x,y) = \int_1^{\e^{-1}} \dj Q_t (x,y) \frac{dt}{\sqrt{t}},
\ \ \ \ \widetilde{\Q_j^\varepsilon}(x,y)= \int_\varepsilon^{1}
\dj \wt{Q_t} (x,y) \frac{dt}{\sqrt{t}},$$
$$\W_j^\e (x,y) = \int_1^{\e^{-1}} \dj W_t (x,y)
\frac{dt}{\sqrt{t}}, \ \ \ \ \ \wt{\W_j^\e} (x,y)= \int_\e^1 \dj
\widetilde{W_t} (x,y) \frac{dt}{\sqrt{t}}.$$
 We shall prove that
$\Q_j^\e$,  $\W_j^\e$, and $\wt{\W_j^\e}$ converge in the norm
operator topology on  $L^1(\mathbb R^d)$, while $\wt{\Q_j^\e}$
converges   strongly on $L^1(\mathbb R^d)$ as $\e$ tends to 0.

\begin{lema} \label{lem1}
The operators $\Q_j^{\e}$ converge in the norm operator topology
on $L^1(\mathbb R^d)$ as $\e\to 0$.
\end{lema}
\begin{proof} There exists $\phi \in \SSS(\R^d)$, $\phi \geq 0$,
such that
 \begin{equation}\label{es1}
\left| \dj P_t(x-z) \right| \leq t^{-1\slash 2} \phi_t(x-z).
 \end{equation}
 On
the other hand, by \eqref{F-Kac}, $K_s(z,y) \leq C s^{-d\slash 2}$.
Hence, for  $0<\varepsilon_2<\varepsilon_1<1$, we have
\begin{equation}\begin{split}
 \int\left| \Q_j^{\e_1}(x,y)-\Q_j^{\e_2}(x,y)\right|dx &
  \leq C \int \int_{\e_1^{-1}}^{\e_2^{-1}} \int_t^\8 \int
t^{-1\slash 2} \phi_t(x-z) V(z) s^{-d\slash 2} dz ds
\frac{dt}{\sqrt{t}} dx \cr &\leq C \int_{\e_1^{-1}}^{\e_2^{-1}}
t^{-d\slash 2} dt \cdot \|V\|_{\LL(\mathbb R^d)},
 \end{split}\end{equation}
 which tents to zero uniformly with respect to $y$ as $\e_1, \
 \e_2\to 0$.

\end{proof}

\begin{lema} \label{lem2}
The operators $\W_j^\e$  converge in the norm operator topology
 on $L^1(\mathbb R^d)$ as $\e\to 0$.
\end{lema}
\begin{proof} The proof borrows ideas from \cite{DZ}.
 Let $0<\varepsilon_2<\varepsilon_1<1$. Then
\begin{align*}
\int & \left|  \W_j^{\e_1} (x,y) -\W_j^{\e_2}(x,y)\right|dx  \\
&\leq \int_{\mathbb R^d} \int_{\e_1^{-1}}^{\e_2^{-1}}
\int_0^{t^{8\slash 9}}\int_{\mathbb R^d}\left| \dj
\Big(P_t(x-z)-P_{t-s}(x-z)\Big)\right| V(z)K_s(z,y)\, dz \, ds
\frac{dt}{\sqrt{t}}\, dx\cr &\ \ + \int_{\mathbb R^d}
\int_{\e_1^{-1}}^{\e_2^{-1}} \int_{t^{8\slash 9}}^t \int_{\mathbb
R^d} \left|\dj \Big(P_t(x-z)-P_{t-s}(x-z)\Big)\right|
V(z)K_s(z,y)\, dz \, ds \frac{dt}{\sqrt{t}} \, dx \cr
&=\W'(y)+\W''(y).
\end{align*}
Observe that there exists $\phi \in \SSS (\R)$, $\phi \geq 0$,
such that for $0<s<t^{8\slash 9}$,
\begin{equation*}
\Big|\dj \Big(P_t(x-z)-P_{t-s} (x-z)\Big)\Big| \leq s \, t^{-3\slash 2}
\phi_t (x-z).
\end{equation*}
Therefore
\begin{align}\label{es3}
\W'(y) &\leq \int \int_{\e_1^{-1}}^{\e_2^{-1}}  \int_0^{t^{8\slash
9}} \int s \, t^{-2} \phi_t(x-z) V(z) K_s(z,y) dz ds dt dx \cr
&\leq \int_{\e_1^{-1}}^{\e_2^{-1}} t^{-10\slash 9} dt \cdot \int
V(z) |z-y|^{2-d} dz \leq C\e_1^{1\slash 9}
\end{align}
uniformly in $y$.  The last inequality is a simple consequence of
the H\"older inequality and the assumption $p > d\slash 2$.

For $t^{8\slash 9} < s <t$ we have $K_s(z,y)\leq C \, s^{-d\slash
2} \leq C \, t^{-4d \slash 9}$. Using (\ref{es1}) we get
\begin{equation}\begin{split}\label{es4} \W''(y)  &\leq C
\int \int_{\e_1^{-1}}^{\e_2^{-1}}  \int_{t^{8\slash 9}}^t \int
\left(t^{-1\slash 2} \phi_t(x-z)+(t-s)^{-1\slash 2}
\phi_{t-s}(x-z)\right) \cr &\ \ \times V(z) t^{-4d \slash 9} dz ds
\frac{dt}{\sqrt{t}} dx\cr &\leq C\|V\|_{\LL(\mathbb R^d)}
\int_{\e_1^{-1}}^{\e_2^{-1}} t^{-4d \slash 9} dt  +
C\|V\|_{\LL(\mathbb R^d)} \int_{\e_1^{-1}}^{\e_2^{-1}} t^{-4d
\slash 9 -1\slash 2} \int_0^t (t-s)^{-1\slash 2} ds \, dt \cr
 & \leq
C\e_1^{4d\slash 9-1}
\end{split}\end{equation}
uniformly in $y$. Now the lemma follows from (\ref{es3}) and
(\ref{es4}).
\end{proof}

\begin{lema} \label{lem4}
There exists a limit of the operators  $\wt{\W_j^\e}$ in the norm
operator topology on $L^1(\mathbb R^d) $ as $\e\to0$.
\end{lema}

\begin{proof}
 Let $0<\e_1<\e_2<1$. Applying  (\ref{es1}) we obtain
\begin{equation}\begin{split}\label{le4}
\int_{\mathbb R^d}
|\widetilde{\W_j^{\e_2}}(x,y)&-\widetilde{\W_j^{\e_2}}(x,y)|\, dx
\\
&\leq \int_{\mathbb R^d} \int_{\e_1}^{\e_2} \int_0^{t}
\int_{\mathbb R^d} (t-s)^{-1\slash 2} \phi_{t-s}(x-z) V (z) K_s
(z,y) \, dz ds \frac{dt}{\sqrt{t}}\, dx \cr &= \int_{\mathbb R^d}
\int_{\e_1}^{\e_2} \int_0^{t\slash 2}\int_{\mathbb R^d}... +
\int_{\mathbb R^d} \int_{\e_1}^{\e_2} \int_{t\slash
2}^t\int_{\mathbb R^d}... =
\widetilde{\W'}(y)+\widetilde{\W''}(y).
\end{split}\end{equation}

If $0 < s < t\slash 2$, then, of course, $(t-s)^{-1\slash 2} \leq
C t^{-1\slash 2}$. Note that
$$\int_0^t K_s(z,y) ds \leq C
|z-y|^{2-d} \exp\left(-\frac{|z-y|^2}{8t}\right) \leq t \psi_t
(z-y)$$
 for some $\psi \in L^{p'} (\R ^d)$, $\psi \geq 0$ ($p'$
denotes the H\"older conjugate exponent to $p$). Hence
\begin{align} \label{lem4a}
 \widetilde{\W'}(y) &\leq C\int_{\e_1}^{\e_2} \int_0^{t\slash
2} \int t^{-1} V(z) K_s(z,y) dz ds dt \leq C\int_{\e_1}^{\e_2}
\int V(z) \psi_t (z-y) dz dt\cr &\leq C\int_{\e_1}^{\e_2} \|V\|_p
\|\psi_t\|_{p'} dt \leq C \int_{\e_1}^{\e_2} t^{-d\slash 2p}
\|\psi_1\|_{p'} dt \leq C\e_2^{1-d\slash 2p}
\end{align}
uniformly in $y$.

If $t\slash 2 \leq s \leq t$, then there exists $\varphi \in
\SSS(\R^d)$, $\varphi \geq 0$, such that $K_s(z,y) \leq
\varphi_t(z-y)$. Therefore
\begin{align}\label{lem4b}
\widetilde{W''}(y)
 &\leq \int_{\e_1}^{\e_2} \int_{t\slash 2}^t
  \int (t-s)^{-1\slash 2} V(z)K_s(z-y) dz ds \frac{dt}{\sqrt{t}} \cr
&\leq C \int_{\e_1}^{\e_2}  \int \int_{0}^{t\slash 2}
(st)^{-1\slash 2} V(z)\varphi_t(z-y) ds dz dt \leq
C\int_{\e_1}^{\e_2} \|V\|_p \|\varphi_t\|_{p'} dt \leq
C\e_2^{1-d\slash 2p}
\end{align}
uniformly in $y$. Now the  lemma is a consequence  \eqref{lem4a}
-- \eqref{lem4b}.
\end{proof}

\begin{lema}\label{lem3}
Assume that $f \in \LL(\mathbb R^d)$. Then the limit $F = \lim_{\e
\to 0} \wt{\Q^\e_j} f$ exists in the $L^1(\mathbb R^d)$-norm.
Moreover, $\|F\|_{\LL(\mathbb R^d)} \leq C \|f\|_{\LL(\mathbb
R^d)}$ with $C$ independent of $f$.
\end{lema}
\begin{proof}
Of course, for any fixed $y \in \R ^d$, the function $ z \mapsto
U(z,y) =  V(z)\Gamma (z,y)$ is supported in the unit ball and
$\|U(z,y)\|_{L^r(dz)} \leq C_r$  for fixed $r\in \left[1,
\frac{dp}{dp+d-2p}\right)$ with $C_r$ independent of $y$. The last
statement follows from  \eqref{F-Kac} and the H\"older inequality.
Let
 $$ H_j^\e(x,z)=\int_\e^1\frac{\partial}{\partial
 x_j}P_t(x-z)\frac{dt}{\sqrt{t}}, \ \ H_j^\e g(x)=\int_{\mathbb
 R^d} H^\e_j(x,z)g(z)\, dz, \ \ H_j^*g(x)=\sup_{0<\e <1}| H_j^\e g(x)|.$$
 It follows from the theory of singular integral convolution
 operators (see, e.g., \cite[Chapter 4]{Duo}) that  for
 $1<r<\infty$ there exists $C_r$ such that
 \begin{equation}\label{eq_max}
 \| H_j^*g\|_{L^r(\mathbb R^d)} \leq C_r\| g\|_{L^r(\mathbb R^d)}\
 \ \ \text{for} \ g\in L^r(\mathbb R^d)
 \end{equation}
 and $\lim_{\e\to 0} H_j^\e g(x)=H_j g(x)$ a.e. and in
 $L^r(\mathbb R^d)$-norm.

 Note that
$\widetilde{\Q_j^\e} (x,y)  =  H^\e_j U(\, \cdot\, ,y)(x)$. Thus
there exists a function $\wt{\Q_j}(x,y)$ such that
 $\lim_{\e\to 0} \wt{\Q_j^\e}(x,y)=\wt{\Q_j}(x,y)$ a.e. and
 \begin{equation}\label{max4}
 \sup_{y}\int_{\mathbb R^d} \sup_{0<\e<1} \left|\wt{\Q^{\e}_j}(x,y)\right|^r\,
 dx\leq C'_r \ \ \text{for} \ 1<r<\frac{dp}{dp+d-2p}.
 \end{equation}
 Since $|H^{\e}_j(x,z)|\leq C_N|x-z|^{-N}$ for $|x-z|>1$,
 \begin{equation}\label{max5}
 |\wt{\Q_j^{\e}}(x,y)|=\left| \int_{|z|\leq 1} H^{\e}_j(x,z)U(z,y)\,
 dz\right|\leq C_N|x|^{-N} \ \ \text{for } |x|>2.
 \end{equation}
 The H\"older inequality combined with (\ref{max4}) and
 (\ref{max5}) implies
  \begin{equation}\label{max6}
 \sup_{y}\int_{\mathbb R^d} \sup_{0<\e<1} \left|\wt{\Q^{\e}_j}(x,y)\right|\,
 dx\leq C,
 \end{equation}
 \begin{equation}\label{limit2}
 \lim_{\e\to 0} \int_{\mathbb R^d}
 \left|\wt{\Q_j^{\e}}(x,y)-\wt{\Q_j}(x,y)\right|\, dx =0 \ \ \text{for every }
 y.
 \end{equation}
 Now the lemma could be easily concluded from (\ref{max6}),
 (\ref{limit2}), and  Lebesgue's dominated convergence theorem.
\end{proof}

%%%%%%%%%%%%%%%%%%%%%%%%%%%%%%%%%%%%%%%%%%%%%%%%%%%%%%%%%%%%%%%%%%%%%%%%%%%%%%%
\section{Proof of the main theorem}

Recall that $(I-VL^{-1})$ is an isomorphism  in $L^1(\mathbb R^d)$.
Consider $f\in \LL(\mathbb R^d)$. Using (\ref{epsilon}) and lemmas
\ref{lem1}, \ref{lem2}, \ref{lem4}, and \ref{lem3}  we get that
$R_jf$ belongs to $L^1(\mathbb R^d)$ if and only if $\mathcal
R_j(I-VL^{-1})f \in L^1(\mathbb R^d)$. Moreover,
\begin{equation}\nonumber
 \| f\|_{L^1(\mathbb R^d)} +\sum_{j=1}^d\| R_j f\|_{L^1(\mathbb
 R^d)} \sim \| (I-VL^{-1})f\|_{L^1(\mathbb R^d)}
 +\sum_{j=1}^d\| \mathcal R_j (I-VL^{-1})f\|_{L^1(\mathbb
 R^d)}.
\end{equation}
Applying the characterization of the classical Hardy space
$H^1(\mathbb R^d)$  by means of the Riesz transforms $\mathcal R_j$
(see (\ref{classical})) and (\ref{homeo}) we obtain the theorem.
%%%%%%%%%%%%%%%%%%%%%%%%%%%   References   %%%%%%%%%%%%%%%%%%%%%%%%%%%%%

\end{document}